\documentstyle[12pt]{article}
\newcommand{\sect}[1]{\section{#1}\setcounter{equation}{0}}

\font\mbn=msbm10 scaled \magstep1
\font\mbs=msbm7 scaled \magstep1
\font\mbss=msbm5 scaled \magstep1
\newfam\mbff
\textfont\mbff=\mbn
\scriptfont\mbff=\mbs
\scriptscriptfont\mbff=\mbss\def\mbf{\fam\mbff}
\def\Re{{\mbf R}}

\def\H{{\mbf H}}
\newtheorem{Th}{Theorem}[section]
\newtheorem{Lm}[Th]{Lemma}

\newtheorem{D}[Th]{Definition}
\newtheorem{Proposition}[Th]{Proposition}

\newtheorem{Problem}[Th]{Problem}
\author{Alexander Brudnyi\thanks{Research supported in part by NSERC.
}\\
Department of Mathematics and Statistics\\
University of Calgary, Calgary,
Canada\\
e-mail: albru@math.ucalgary.ca\\
\\
Yuri Brudnyi\\
Department of Mathematics\\
Technion, Haifa,
Israel\\
e-mail: ybrudnyi@math.technion.ac.il}
\title{Extension of Lipschitz Functions Defined on Metric Subspaces of Homogeneous Type}
\date{}
\begin{document}
\maketitle
\begin{abstract}
{If a metric subspace $M^{o}$ of an arbitrary metric space $M$ carries a doubling measure $\mu$, then there is a simultaneous linear extension of all Lipschitz functions on $M^{o}$ ranged in a Banach space to those on $M$. Moreover, the norm of this linear operator is controlled by logarithm of the doubling constant of $\mu$.}
\end{abstract}
\sect{Formulation of the Main Result}
Let $(M,d)$ be a metric space and $X$ be a Banach space.
The space $Lip(M,X)$ consists of all $X$-valued Lipschitz functions on $M$.
The Lipschitz constant
\begin{equation}\label{e1}
L(f):=\sup_{m\neq m'}\left\{\frac{||f(m)-f(m')||}{d(m,m')}\ :\ m, m'\in M\ \right\}
\end{equation}
of a function $f$ from this space is therefore finite and the function $f\mapsto L(f)$ is a Banach seminorm on $Lip(M,X)$.

Let $M^{o}$ be a metric subspace of $M$, i.e., $M^{o}\subset M$ is a metric space endowed with the induced metric $d|_{M^{o}\times M^{o}}$.\\
{\underline{\em Convention.}} We mark all objects related to the subspace $M^{o}$ by the upper ``$o$``.

A linear operator $E: Lip(M^{o},X)\to Lip(M,X)$ is called a
{\em simultaneous extension} if for all $f\in Lip(M^{o},X)$
$$
Ef|_{M^{o}}=f
$$
and, moreover, the norm
$$
||E||:=\sup\left\{\frac{L(Ef)}{L(f)}\ :\ f\in Lip(M^{o},X)\right\}
$$
is finite.

To formulate the main result we also need
\begin{D}\label{de1}
A Borel measure $\mu$ on a metric space $(M,d)$ is said to be doubling if the $\mu$-measure of every open ball
$$
B_{R}(m):=\{m'\in M\ :\ d(m,m')<R\}
$$
is strictly positive and finite and the doubling constant
\begin{equation}\label{eq1}
D(\mu):=\sup\left\{\frac{\mu(B_{2R}(m))}{\mu(B_{R}(m))}\ :\ m\in M, \ R>0\right\}
\end{equation}
is finite.

A metric space carrying a fixed doubling measure is called of homogeneous type.
\end{D}

Our main result is
\begin{Th}\label{t1}
Let $M^{o}$ be a metric subspace of an arbitrary metric space $(M,d)$. Assume that $(M^{o},d^{o})$ is of homogeneous type and $\mu^{o}$ is the corresponding doubling measure. Then there exists a simultaneous extension
$E:Lip(M^{o},X)\to Lip(M,X)$ satisfying
\begin{equation}\label{eq2}
||E||\leq c(\log_{2} D(\mu^{o})+1)
\end{equation}
with some numerical constant $c>1$.
\end{Th}

Let us discuss relations of this theorem to some known results. First, a similar result holds for an arbitrary subspace $M^{o}$ provided that the ambient space $M$ is of {\em pointwise homogeneous type}, see [BB1, Theorem 2.21] and [BB2, Theorem 1.14]. The class of metric spaces of pointwise homogeneous type contains, in particular, all metric spaces of homogeneous type, Riemannian manifolds $M_{\omega}\cong\Re^{n}\times\Re_{+}$ with the path metric defined by the Riemannian metric
$$
ds^{2}:=\omega(x_{n+1})(dx_{1}^{2}+\dots+dx_{n+1}^{2}),\ \ \
(x_{1},\dots, x_{n},x_{n+1})\in\Re^{n}\times\Re_{+}\ \!,
$$
where $\omega:\Re_{+}\to\Re_{+}$ is a continuous nonincreasing function
(e.g., the hyperbolic spaces $\H^{n}$ are in this class), and finite direct products of these objects.

The following problem is of a considerable interest.

\begin{Problem}\label{p1}
Is it true that Theorem \ref{t1} is valid for $M^{o}\ \!(\subset \!M)$ isometric to a subspace of a metric space $(\widehat M, \widehat d)$ of pointwise homogeneous type with $||E||\leq c(\widehat M)$?
\end{Problem}
(Here $c(\widehat M)$ depends on some characteristics of
$\widehat M$ only.)

It is proved in [BB2] that as such $M^{o}$ one can take, e.g., finite direct products of Gromov hyperbolic spaces of bounded geometry and that 
the answer in Problem \ref{p1} is positive in this case.
 
Second, as a consequence of Theorem \ref{t1} we obtain a deep extension result due to Lee and Naor, see [LN, Theorem 1.6]. The latter asserts that a simultaneous extension $E:Lip(M^{o},X)\to Lip(M,X)$ exists whenever the subspace $(M^{o},d^{o})$ of $(M,d)$ has the finite doubling constant $\delta(M^{o})$ and, moreover,
\begin{equation}\label{eq3}
||E||\leq c\log_{2}\delta(M^{o})
\end{equation}
with some numerical constant $c>1$.

Let us recall that the doubling constant $\delta(M)$ of a metric space $(M,d)$ is the infimum of integers $N$ such that every closed ball of $M$ of radius $R$ can be covered by $N$ closed balls of radius $R/2$. The space $M$ is said to be {\em doubling} if $\delta(M)<\infty$.

To derive the Lee-Naor theorem from our main result we first note that without loss of generality one may assume that $(M^{o},d^{o})$ is complete. By the Koniagin-Vol'berg theorem [KV] (see also [LS]) a complete doubling space $M$ carries a doubling measure $\mu$ such that
\begin{equation}\label{eq4}
\log_{2}D(\mu)\leq c\log_{2}\delta(M)
\end{equation}
where $c\geq 1$ is a numerical constant. Together with  (\ref{eq2}) this implies the Naor-Lee result.

On the other hand, it was noted in [CW] that if $M$ carries a doubling measure $\mu$, then this space is doubling and
\begin{equation}\label{eq5}
\log_{2}\delta(M)\leq c\log_{2}D(\mu)
\end{equation}
with some numerical constant $c>1$.
Hence, Theorem \ref{t1} is, in turn, a consequence of (\ref{eq5}) and the Lee-Naor theorem. However, the rather elaborated proof of the latter result is nonconstructive (it exploits an appropriate stochastic metric decomposition of $M\setminus M^{o}$). In contrast, our proof is constructive and is based on a simple average procedure. Therefore our proof can be also seen as a streamlining constructive method of the proof of the Lee-Naor theorem.
\sect{\hspace*{-1em}. Proof of Theorem \ref{t1}.}
We begin with the following remark reducing the required result to a special case.

Let $M$ and $M^{o}$ be isometric to subspaces of a new metric space $\widehat M$ and its subspace $\widehat M^{o}$, respectively. Assume that there exists a simultaneous extension $\widehat E: Lip(\widehat M^{o},X)\to Lip(\widehat M,X)$. Then, after identification of $M^{o}$ and $M$ with the corresponding isometric subspaces of $\widehat M$, the operator $\widehat E$ gives rise to a simultaneous extension $E:Lip(M^{o},X)\to Lip(M,X)$ satisfying
\begin{equation}\label{eq6}
||E||\leq ||\widehat E||.
\end{equation}
If, in addition, $||\widehat E||$ is bounded by the right-hand side of (\ref{eq2}), then the desired result immediately follows.

We choose as the above pair $\widehat M^{o}\subset\widehat M$ metric spaces
denoted by $M_{N}^{o}$ and $M_{N}$ where $N\geq 1$ is a fixed integer and defined as follows.

The underlying sets of these spaces are
\begin{equation}\label{eq7}
M_{N}:=M\times\Re^{N},\ \ \ M_{N}^{o}:=M^{o}\times\Re^{N};
\end{equation}
a metric $d_{N}$ on $M_{N}$ is given by
\begin{equation}\label{eq8}
d_{N}((m,x),(m',x')):=d(m,m')+|x-x'|_{1}
\end{equation}
where $m,m'\in M$ and $x,x'\in\Re^{N}$, and
$|x|_{1}:=\sum_{i=1}^{N}|x_{i}|$ is the $l_{1}^{N}$-metric of $x\in\Re^{N}$.
Further, $d_{N}^{o}$ denotes the metric on $M_{N}^{o}$ induced by $d_{N}$.

Finally, we define a Borel measure $\mu_{N}^{o}$ on $M_{N}^{o}$ as the tensor product of the measure $\mu^{o}$ and the Lebesgue measure $\lambda_{N}$ on $\Re^{N}$:
\begin{equation}\label{eq9}
\mu_{N}^{o}:=\mu^{o}\otimes\lambda_{N}.
\end{equation}
We extend this measure to the $\sigma$-algebra consisting of subsets $S\subset M_{N}$ such that $S\cap M_{N}^{o}$ is a Borel subset of $M_{N}^{o}$. Namely, we set for these $S$
$$
\overline\mu_{N}(S):=\mu_{N}^{o}(S\cap M_{N}^{o}).
$$
It is important for the subsequent part of the proof that every open ball $B_{R}((m,x))\subset M_{N}$ belongs to this $\sigma$-algebra. In fact, its intersection with $M_{N}^{o}$ is a Borel subset of this space, since the function $(m',x')\mapsto d_{N}((m,x), (m',x'))$ is continuous on $M_{N}^{o}$. Hence,
\begin{equation}\label{eq10}
\overline\mu_{N}(B_{R}((m,x)))=\mu_{N}^{o}(B_{R}((m,x))\cap M_{N}^{o}).
\end{equation}

{\underline{\em Auxiliary results}}

The measure $\mu_{N}^{o}$ is clearly doubling. Therefore its {\em dilation function} given for $l\geq 1$ by
$$
D_{N}^{o}(l):=\sup\left\{\frac{\mu_{N}^{o}(B_{lR}^{o}(\widehat m))}{\mu_{N}^{o}(B_{R}^{o}(\widehat m))}\ :\ \widehat m\in M_{N}^{o}\ {\rm and}\ R>0\right\}
$$
is finite.\\
Hereafter we denote by $\widehat m$ the pair $(m,x)$ with $m\in M$ and $x\in\Re^{N}$, and by $B_{R}^{o}(\widehat m)$ the open ball in $M_{N}^{o}$ centered at $\widehat m\in M_{N}^{o}$ and of radius $R$. The open ball $B_{R}(\widehat m)$ of $M_{N}$ relates to that by
$$
B_{R}^{o}(\widehat m)=B_{R}(\widehat m)\cap M_{N}^{o}
$$
provided $\widehat m\in M_{N}^{o}$.

In [BB1] the value $D_{N}^{o}(1+1/N)$ is proved to be bounded by some numerical constant for all sufficiently large $N$. In the argument presented below we require a similar estimate for a (modified) dilation function $D_{N}$ for the extended measure $\overline\mu_{N}$. This is given for $l\geq 1$ by
\begin{equation}\label{eq11}
D_{N}(l):=\sup\left\{\frac{\overline\mu_{N}(B_{lR}(\widehat m))}{\overline\mu_{N}(B_{R}(\widehat m))}\right\}
\end{equation}
where the supremum is taken over all $R$ satisfying
\begin{equation}\label{eq12}
R>4d(\widehat m,M_{N}^{o}):=4\inf\{d_{N}(\widehat m,\widehat m')\ :\ \widehat m'\in M_{N}^{o}\}
\end{equation}
and then over all $\widehat m\in M_{N}$.

Due to (\ref{eq10}) and (\ref{eq12}) the denominator in (\ref{eq11}) is not zero and $D_{N}(l)$ is well defined. 

Comparison of the above dilation functions shows that $D_{N}^{o}(l)\leq D_{N}(l)$. Nevertheless, the converse is also true for $l$ close to $1$.
\begin{Lm}\label{l128}
Assume that $N$ and the doubling constant $D:=D(\mu^{o})$, see (\ref{eq1}), are related by
\begin{equation}\label{eq13}
N\geq [3\log_{2} D]+5\ .
\end{equation}
Then the following is true:
$$D_{N}(1+1/N)\leq \frac{6}{5}e^{4}\ .
$$
\end{Lm}
{\bf Proof.}
In accordance with the definition of $D_{N}$, see (\ref{eq11}), we must estimate the function
\begin{equation}\label{eq15}
\frac{\overline\mu_{N}(B_{R_{N}}(\widehat m))}{\overline\mu_{N}(B_{R}(\widehat m))}\ \ \ {\rm where}\ \ \ R_{N}:=\left(1+\frac{1}{N}\right)R.
\end{equation}
Since the points $\widehat m'$ of the ball $B_{R_{N}}(\widehat m)$ of $M_{N}$ satisfy the inequality
$$
d(m,m')+|x-x'|_{1}<R_{N},
$$
the Fubini theorem and (\ref{eq10}) yield
\begin{equation}\label{eq16}
\overline\mu_{N}(B_{R_{N}}(\widehat m))=\gamma_{N}\int_{M^{o}\cap B_{R_{N}}(m)}
(R_{N}-d(m,m'))^{N}d\mu^{o}(m')\ \!;
\end{equation}
here $\gamma_{N}$ is the volume of the 
unit $l_{1}^{N}$-ball. 

We must estimate the integral in (\ref{eq16}) from above under the condition
\begin{equation}\label{eq17}
d_{N}(\widehat m,M_{N}^{o})<R/4.
\end{equation}
To this end 
split the integral into one over $B_{3R/4}(m)\cap M^{o}$ and one
over the remaining part $(B_{R_{N}}(m)\setminus B_{3R/4}(m))\cap M^{o}$. Denote these
integrals by $I_{1}$ and $I_{2}$. For $I_{2}$ we get 
\begin{equation}\label{eq18}
I_{2}\leq \gamma_{N}(R_{N}-3R/4)^{N}\mu^{o}(B_{R_{N}}(m)\cap M^{o}).
\end{equation}
Further, from (\ref{eq17}) we clearly have
$$
d(m,M^{o})<R/4.
$$
Pick a point $\widetilde m\in M^{o}$ so that
$$
d(m,M^{o})\leq d(m,\widetilde m)< R/4.
$$
Then we have the following embeddings
$$
B_{R_{N}/4}^{o}(\widetilde m)\subset B_{R_{N}/2}(m)\cap M^{o}\subset B_{R_{N}}(m)\cap M^{o}\subset B_{5R_{N}/4}^{o}(\widetilde m).
$$
Applying the doubling inequality for the measure $\mu^{o}$, see (\ref{eq1}), we then obtain
$$
\mu^{o}(B_{R_{N}}(m)\cap M^{o})\leq D^{3}\mu^{o}(B_{R_{N}/2}(m)\cap M^{o})\ .
$$
Moreover, due to (\ref{eq13}) 
$$
D^{3}<2^{[3\log_{2}D]+1}\leq 2^{N-4}.
$$
Combining the last two inequalities with (\ref{eq18}) we have
\begin{equation}\label{eq19}
I_{2}\leq\gamma_{N}2^{-N-4}\left(1+\frac{4}{N}\right)^{N}R^{N}
\mu^{o}(B_{R_{N}/2}(m)\cap M^{o})\ .
\end{equation}
To estimate the integral $I_{1}$ we rewrite its integrand as follows:
$$
(R_{N}-d(m,m'))^{N}=\left(1+\frac{1}{N}\right)^{N}(R-d(m,m'))^{N}\left(1+
\frac{d(m,m')}{(N+1)(R-d(m,m'))}\right)^{N}.
$$
Since $m'\in B_{3/4R}(m)$, the last factor is 
at most $\left(1+\frac{3R/4}{(N+1)R/4}\right)^{N}=\left(1+\frac{3}{N+1}\right)^{N}$. This yields
$$
\begin{array}{c}
\displaystyle
I_{1}\leq\gamma_{N}
\left(1+\frac{1}{N}\right)^{N}\left(1+\frac{3}{N+1}\right)^{N}
\int_{B_{3R/4}(m)\cap M^{o}}(R-d(m,m'))^{N}d\mu^{o}(m')\leq\\
\\
\displaystyle
 e^{4}\overline\mu_{N}(B_{R}(\widehat m)).
\end{array}
$$
Hence for the part of fraction (\ref{eq15}) related to $I_{1}$ we have
\begin{equation}\label{eq20}
\widetilde I_{1}:=\frac{I_{1}}{\overline\mu_{N}(B_{R}(\widehat m))}\leq e^{4}.
\end{equation}
To estimate the remaining part $\widetilde I_{2}:=\frac{I_{2}}{\overline\mu_{N}(B_{R}(\widehat m))}$ we note that its denominator is greater than
$$
\gamma_{N}\int_{M^{o}\cap B_{R_{N}/2}(m)}(R-d(m,m'))^{N}d\mu^{o}(m').
$$
Since here 
$d(m,m')\leq R_{N}/2$, this, in turn, is bounded from below by
$$
\gamma_{N}2^{-N}\left(1-\frac{1}{N}\right)^{N}R^{N}
\mu^{o}(B_{R_{N}/2}(m)\cap M^{o}). 
$$
Combining this with (\ref{eq19}) and noting that $N\geq 5$ we get
$$
\widetilde I_{2}\leq 2^{-4}\left(1-\frac{1}{N}\right)^{-N}
\left(1+\frac{4}{N}\right)^{N}<\frac{1}{5}e^{4}.
$$
Hence the fraction (\ref{eq15}) is bounded by $\widetilde I_{1}+\widetilde I_{2}\leq\frac{6}{5}e^{4}$, see (\ref{eq20}), and this immediately implies the required estimate of
$D_{N}(1+1/N)$.\ \ \ \ \ $\Box$

In the next lemma we estimate $\overline\mu_{N}$-measure of the spherical layer \penalty-10000
$B_{R_{2}}(\widehat m)-B_{R_{1}}(\widehat m)$, $R_{2}\geq R_{1}$, by a
kind of a surface measure. For its formulation we set
\begin{equation}\label{eq21}
A_{N}:=\frac{12}{5}e^{4}N.
\end{equation}
\begin{Lm}\label{l3.5}
Assume that
$$
N\geq [3\log_{2}D]+6.
$$
Then for all $\widehat m\in M_{N}$ and $R_{1},R_{2}>0$ satisfying 
$$
R_{2}\geq\max\{R_{1}, 8d_{N}(\widehat m,M_{N}^{o})\} 
$$
the following is true
$$
\overline\mu_{N}(B_{R_{2}}(\widehat m)\setminus B_{R_{1}}(\widehat m))\leq A_{N}
\frac{\overline\mu_{N}(B_{R_{2}}(\widehat m))}{R_{2}}(R_{2}-R_{1}).
$$
\end{Lm}
{\bf Proof.} By definition $M_{N}=M_{N-1}\times\Re$ and 
$\overline\mu_{N}=\overline\mu_{N-1}\otimes\lambda_{1}$. Then by the
Fubini theorem we have for $R_{1}\leq R_{2}$ with $\widehat m=(\widetilde m, t)$
$$
\overline\mu_{N}(B_{R_{2}}(\widehat m))-\overline\mu_{N}(B_{R_{1}}
(\widehat m))
=2\int_{R_{1}}^{R_{2}}\overline\mu_{N-1}(B_{s}(\widetilde m))ds\leq
\frac{2R_{2}\overline\mu_{N-1}(B_{R_{2}}(\widetilde m))}{R_{2}}(R_{2}-R_{1}).
$$
We claim that for arbitrary $l>1$ and $R\geq 8d_{N}(\widehat m, M_{N}^{o}):=8d_{N-1}(\widetilde m, M_{N-1}^{o})$
\begin{equation}\label{e1242}
R\overline\mu_{N-1}(B_{R}(\widetilde m))\leq\frac{lD_{N-1}(l)}{l-1}
\overline\mu_{N}(B_{R}(\widehat m))\ .
\end{equation}
Together with the previous inequality this will yield
$$
\overline\mu_{N}(B_{R_{2}}(\widehat m))-
\overline\mu_{N}(B_{R_{1}}(\widehat m))\leq\frac{2lD_{N-1}(l)}{l-1}\cdot
\frac{\overline\mu_{N}(B_{R_{2}}(\widehat m))}{R_{2}}(R_{2}-R_{1}).
$$
Finally choose here $l=1+\frac{1}{N-1}$ and use Lemma \ref{l128}. This will give the required inequality.

Hence, it remains to establish (\ref{e1242}). By the definition of
$D_{N-1}(l)$ we have for $l>1$ using the previous lemma
$$
\begin{array}{c}
\displaystyle
\overline\mu_{N}(B_{lR}(\widehat m))=2l\int_{0}^{R}
\overline\mu_{N-1}(B_{ls}(\widetilde m))ds\leq 4l\int_{R/2}^{R}\overline\mu_{N-1}(B_{ls}(\widetilde m))ds\leq\\
\\
\displaystyle 4lD_{N-1}(l)\int_{R/2}^{R}\overline\mu_{N-1}(B_{s}(\widetilde m))ds\leq
2lD_{N-1}(l)\overline\mu_{N}
(B_{R}(\widehat m))\ .
\end{array}
$$
On the other hand, replacing $[0,R]$ by $[l^{-1}R,R]$ we also have
$$
\overline\mu_{N}(B_{lR}(\widehat m))\geq
2l\overline\mu_{N-1}(B_{R}(\widetilde m))
(R-l^{-1}R)= 2(l-1)R\overline\mu_{N-1}(B_{R}(\widetilde m))\ .
$$
Combining the last two inequalities we get (\ref{e1242}).\ \ \ \ \ $\Box$
\vspace{3mm}

{\underline{\em Extension operator}}

We define the required simultaneous extension $E:Lip(M_{N}^{o},X)\to Lip(M_{N},X)$ using the standard average operator $Ave$ defined on continuous and locally bounded functions $g:M_{N}^{o}\to X$ by
$$
Ave(g;\widehat m,R):=\frac{1}{\overline\mu_{N}(B_{R}(\widehat m))}\int_{B_{R}(\widehat m)}g\ \!d\overline\mu_{N}.
$$

To be well-defined the domain of integration $B_{R}(\widehat m)\cap M_{N}^{o}$ should be of strictly positive $\overline\mu_{N}$-measure (i.e., $\mu_{N}^{o}$-measure). This condition is fulfilled in the case presented now. Namely, we define the simultaneous extension $E$ on functions $f\in Lip(M_{N}^{o},X)$ by
\begin{equation}\label{e3.13}
(Ef)(\widehat m):=\left\{
\begin{array}{ccc}
f(\widehat m)&{\rm if}&\widehat m\in M_{N}^{o}\\
\\
Ave(f;m, R(\widehat m))&{\rm if}&\widehat m\notin M_{N}^{o}
\end{array}
\right.
\end{equation}
where we set
$$
R(\widehat m):=8d_{N}(\widehat m,M_{N}^{o}).
$$
The required estimate of $||E||$ is presented below. To formulate the result we set
\begin{equation}\label{e3.29}
K_{N}(l):=A_{N}D_{N}(l)(4l+1)
\end{equation}
where the first of two factors are defined by (\ref{eq21}) and (\ref{eq11}).
\begin{Proposition}\label{p3.6}
The following inequality
$$
||E||\leq 20A_{N}+\max\left(\frac{4l+1}{2(l-1)},K_{N}(l)\right)
$$
is true provided $l:=1+1/N$.
\end{Proposition}

Before we begin the proof let us derive from here the desired result. 
Namely, choose
$$
N:=[3\log_{2} D]+6
$$
and use Lemma \ref{l128} and (\ref{eq21}) to estimate $D_{N}(1+1/N)$ and
$A_{N}$.
Then we get  
$$
||E||\leq C(\log_{2} D+2)
$$
with some numerical constant $C$. This clearly gives (\ref{eq2}).\\
{\bf Proof.}
We have to show that for every $\widehat m_{1},\widehat m_{2}\in  M_{N}$
\begin{equation}\label{e148}
||(Ef)(\widehat m_{1})-(Ef)(\widehat m_{2})||_{X}\leq 
K||f||_{Lip(M_{N},X)}d_{N}(\widehat m_{1},\widehat m_{2})
\end{equation}
where $K$ is the constant in the inequality of 
the proposition.\\
It suffices to consider only two cases:
\begin{itemize}
\item[(a)]
$\widehat m_{1}\in M_{N}^{o}$\ and\ $\widehat m_{2}\not\in M_{N}^{o}$;
\item[(b)]
$\widehat m_{1}, \widehat m_{2}\not\in M_{N}^{o}$.
\end{itemize}
We assume without loss of generality that
\begin{equation}\label{e149}
||f||_{Lip(M_{N}^{o},X)}=1
\end{equation}
and simplify the computations by introducing the following notations:
\begin{equation}\label{e1410}
R_{i}:=d_{N}(\widehat m_{i},M_{N}^{o})\ ,\ 
B_{ij}:=B_{8R_{j}}(\widehat m_{i})\ ,
\ v_{ij}:=\overline\mu_{N}(B_{ij})\ ,\ \ 1\leq i,j\leq 2\ .
\end{equation}
We assume also for definiteness that
\begin{equation}\label{e1411}
0< R_{1}\leq R_{2}\ .
\end{equation}
By the triangle inequality we then have
\begin{equation}\label{e1412}
0\leq R_{2}-R_{1}\leq d_{N}(\widehat m_{1},\widehat m_{2})\ .
\end{equation}
Further, by Lemma \ref{l3.5} the quantities introduced satisfy the following inequality:
\begin{equation}\label{e1414}
v_{i2}-v_{i1}\leq\frac{A_{N}v_{i2}}{R_{2}}(R_{2}-R_{1})\ ,
\end{equation}

Let now $\widehat m^{*}$ be such that $d_{N}(\widehat m_{1}, \widehat m^{*})<2R_{1}$. Set
\begin{equation}\label{e1418}
\widehat f(\widehat m):=f(\widehat m)-
f(\widehat m^{*}).
\end{equation}
From the triangle inequality we then obtain
\begin{equation}\label{e1417}
\max\{||\widehat f(\widehat m)||_{X}\ :\ \widehat m\in B_{i2}\cap M_{N}^{o}\}
\leq 10 R_{2}+(i-1)d_{N}(\widehat m_{1},\widehat m_{2});
\end{equation}
here $i=1,2$.

We now prove (\ref{e148}) for $\widehat m_{1}\in M_{N}^{o}$ and 
$\widehat m_{2}\not\in M_{N}^{o}$. We begin with the evident inequality
$$
||(Ef)(\widehat m_{2})-(Ef)(\widehat m_{1})||_{X}=
\frac{1}{v_{22}}\left|\left|
\int_{B_{22}}\widehat f(\widehat m)d\overline\mu_{N}\right|\right|_{X}
\leq\max_{B_{22}\cap M_{N}^{o}}||\widehat f||_{X},
$$
see (\ref{e1410}) and (\ref{e1418}).
Applying (\ref{e1417}) with $i=2$ we then bound this maximum by
$10R_{2}+d_{N}(\widehat m_{1},\widehat m_{2})$. But 
$\widehat m_{1}\in M_{N}^{o}$ and so
$$
R_{2}=d_{N}(\widehat m_{2},M_{N}^{o})\leq d_{N}(\widehat m_{1},\widehat m_{2});
$$
therefore (\ref{e148}) holds in this case with $K=11$.

The remaining case $\widehat m_{1},\widehat m_{2}\not\in M_{N}^{o}$ requires 
some additional auxiliary results. For their formulations we first write
\begin{equation}\label{e1419}
(Ef)(\widehat m_{1})-(Ef)(\widehat m_{2}):=D_{1}+D_{2}
\end{equation}
where
\begin{equation}\label{e1420}
\begin{array}{c}
D_{1}:=Ave(\widehat f; \widehat m_{1},8R_{1})-
Ave(\widehat f; \widehat m_{1},8R_{2})\\
\\
D_{2}:=Ave(\widehat f; \widehat m_{1},8R_{2})-
Ave(\widehat f; \widehat m_{2},8R_{2})\ ,
\end{array}
\end{equation}
see (\ref{e3.13}) and (\ref{e1418}).
\begin{Lm}\label{l143}
We have
$$
||D_{1}||_{X}\leq 20A_{N}d_{N}(\widehat m_{1},\widehat m_{2})\ .
$$
Recall that $A_{N}$ is the constant defined by (\ref{eq21}).
\end{Lm}
{\bf Proof.}
By (\ref{e1420}), (\ref{e1418}) and (\ref{e1410}),
$$
D_{1}=\frac{1}{v_{11}}\int_{B_{11}}\widehat f d\overline\mu_{N}-
\frac{1}{v_{12}}\int_{B_{12}}\widehat f d\overline\mu_{N}=
\left(\frac{1}{v_{11}}-\frac{1}{v_{12}}\right)\int_{B_{11}}
\widehat f d\overline\mu_{N}-\frac{1}{v_{12}}\int_{B_{12}\setminus B_{11}}
\widehat f d\overline\mu_{N}.
$$
This immediately implies that
$$
||D_{1}||_{X}\leq 2\cdot \frac{v_{12}-v_{11}}{v_{12}}\cdot
\max_{B_{12}\cap M_{N}^{o}}||\widehat f||_{X}\ .
$$
Applying now (\ref{e1414}) and (\ref{e1412}), and then 
(\ref{e1417}) with $i=1$ we get the desired estimate.\ \ \ \ \ $\Box$

To obtain a similar estimate for $D_{2}$ we will use the following
two facts.
\begin{Lm}\label{l144}
Assume that for a given $l>1$
\begin{equation}\label{e1421}
d_{N}(\widehat m_{1},\widehat m_{2})\leq 8(l-1)R_{2}\ .
\end{equation}
Let for definiteness
\begin{equation}\label{eq1421}
v_{22}\leq v_{12}\ .
\end{equation}
Then we have
\begin{equation}\label{e1422}
\overline\mu_{N}(B_{12}\Delta B_{22})\leq A_{N}D_{N}(l)\frac{v_{12}}{4R_{2}}
d_{N}(\widehat m_{1},\widehat m_{2})
\end{equation}
(here $\Delta$ denotes symmetric difference of sets).
\end{Lm}
{\bf Proof.}
Set
$$
R:=8R_{2}+d_{N}(\widehat m_{1},\widehat m_{2})\ .
$$
Then $B_{12}\cup B_{22}\subset B_{R}(\widehat m_{1})\cap 
B_{R}(\widehat m_{2})$, and
\begin{equation}\label{eq1223}
\begin{array}{c}
\displaystyle
\overline\mu_{N}(B_{12}\Delta B_{22})\leq (\overline\mu_{N}(B_{R}(\widehat m_{1}))-\overline\mu_{N}(B_{8R_{2}}(\widehat m_{1})))
+
\\
\\
\displaystyle
(\overline\mu_{N}(B_{R}(\widehat m_{2}))-\overline\mu_{N}(B_{8R_{2}}(\widehat m_{2})))
.
\end{array}
\end{equation}
Estimating the terms on the right-hand side by
Lemma \ref{l3.5} we bound them by 
$$
A_{N}\frac{\overline\mu_{N}(B_{R}(\widehat m_{1}))}{R}\ \! (R-8R_{2})+A_{N}\frac{\overline\mu_{N}(B_{R}(\widehat m_{2}))}{R}(R-8R_{2})\ .
$$
Moreover, $8R_{2}\leq R\leq 8lR_{2}$ and $R-8R_{2}:=
d_{N}(\widehat m_{1},\widehat m_{2})$,
see (\ref{e1421}); taking into account (\ref{eq11}),
(\ref{e1410}) and (\ref{eq1421}) we therefore have
$$
\overline\mu_{N}(B_{12}\Delta B_{22})\leq 
A_{N}D_{N}(l)\frac{v_{12}}{4R_{2}}\ \! 
d_{N}(\widehat m_{1},\widehat m_{2})\ .\ \ \ \ \ \Box
$$
\begin{Lm}\label{l145}
Under the assumptions of the previous lemma we have
\begin{equation}\label{e1423}
v_{12}-v_{22}\leq A_{N}D_{N}(l)\frac{v_{12}}{4R_{2}}
d_{N}(\widehat m_{1},\widehat m_{2})\ .
\end{equation}
\end{Lm}
{\bf Proof.} By (\ref{e1410}) the left-hand side is bounded by
$\overline\mu_{N}(B_{12}\Delta B_{22})$.
 \ \ \ \ \ $\Box$

We now estimate $D_{2}$ from (\ref{e1420}) beginning with
\begin{Lm}\label{l146}
Under the conditions of Lemma \ref{l144} we have
$$
||D_{2}||_{X}\leq K_{N}(l)d_{N}(\widehat m_{1},\widehat m_{2})
$$
where $K_{N}(l):=A_{N}D_{N}(l)(4l+1)$.
\end{Lm}
{\bf Proof.} By the definition of $D_{2}$ and our notation, see
(\ref{e1420}), (\ref{e1418}) and (\ref{e1410}),
$$
\begin{array}{c}
\displaystyle
||D_{2}||_{X}:=\left|\left|\frac{1}{v_{12}}\int_{B_{12}}\widehat f d\overline\mu_{N}-
\frac{1}{v_{22}}\int_{B_{22}}\widehat f d\overline\mu_{N}\right|\right|_{X}\leq
\\
\\
\displaystyle
\frac{1}{v_{12}}\int_{B_{12}\Delta B_{22}}||\widehat f||_{X}\ \! d\overline\mu_{N}
+\left|\frac{1}{v_{12}}-\frac{1}{v_{22}}\right|\int_{B_{22}}
||\widehat f||_{X}
\ \! d\overline\mu_{N}:=J_{1}+J_{2}\ .
\end{array}
$$
By (\ref{e1422}), (\ref{e1421}) and (\ref{e1417})
$$
\begin{array}{c}
\displaystyle J_{1}\leq\frac{1}{v_{12}}\overline\mu_{N}(B_{12}\Delta
B_{22})\sup_{(B_{12}\Delta B_{22})\cap M_{N}^{o}}||\widehat f||_{X}\leq\\
\\
\displaystyle
\frac{A_{N}D_{N}(l)}{4R_{2}}
d_{N}(\widehat m_{1},\widehat m_{2})
(d_{N}(\widehat m_{1},\widehat m_{2})+10R_{2})\leq
A_{N}D_{N}(l)(2l+1/2)d_{N}(\widehat m_{1},\widehat m_{2}).
\end{array}
$$
Also, (\ref{e1423}), (\ref{e1417}) and (\ref{e1421}) yield
$$
J_{2}\leq A_{N}D_{N}(l)(2l+1/2)
d_{N}(\widehat m_{1},\widehat m_{2})\ .
$$
Combining these we get the required estimate.\ \ \ \ \ $\Box$

It remains to consider the case of 
$\widehat m_{1},\widehat m_{2}\in M_{N}$ satisfying
the inequality
$$
d_{N}(\widehat m_{1},\widehat m_{2})>8(l-1)R_{2}
$$
converse to (\ref{e1421}). Now the definition (\ref{e1420}) of
$D_{2}$ and (\ref{e1417}) imply that
$$
||D_{2}||_{X}\leq 2\sup_{(B_{12}\cup B_{22})\cap M_{N}^{o}}||\widehat f||_{X}\leq
2(10R_{2}+d_{N}(\widehat m_{1},\widehat m_{2}))\leq 
\frac{4l+1}{2(l-1)}d_{N}(\widehat m_{1},\widehat m_{2})\ .
$$
Combining this with the inequalities of Lemmas \ref{l143} and
\ref{l146} and equality (\ref{e1419}) we obtain the required
estimate of the Lipschitz norm of the extension operator $E$:
$$
||E||\leq  20A_{N}+ \max\left(\frac{4l+1}{2(l-1)},K_{N}(l)\right)
$$
where $K_{N}(l)$ is the constant in (\ref{e3.29}). \ \ \ \ \ $\Box$


\end{document}